\documentclass[12pt]{article}
\usepackage{amsthm,amssymb}
\usepackage{amsbsy,amsfonts,amsmath}
\newfont{\ssbe}{cmssbx10}

\newcommand{\al}{{\alpha}}

\newcommand{\Aoneoneone}{{A(1,1)^{(1)}}}

\newcommand{\slninione}{{{\rm sl}(2|2)^{(1)}}}

\newcommand{\Uqslninione}{U_q{\rm sl}(2|2)^{(1)}}
\newcommand{\G}{{\frak{g}}}\newcommand{\ovG}{{\overline {\frak{g}}}}

\newcommand{\C}{{\Bbb{C}}}\newcommand{\Z}{{\Bbb{Z}}}
\newcommand{\Ct}{{\Bbb{C}}[t,t^{-1}]}

\newcommand{\Ctms}{{\C^\times}}
\newcommand{\D}{{{\frak{d}}}}

\newcommand{\UqDone}{{{U_q{\D}}^{(1)}}}
\newcommand{\wtUqDone}{{{\widetilde U}_q{\D}^{(1)}}}
\newcommand{\wtU}{{{\widetilde U}}}
\newcommand{\wtUp}{{{\widetilde U}^\prime}}
\newcommand{\wtUzero}{{{\widetilde U}^0}}
\newcommand{\wtUqD}{{{\widetilde U}_q{\D}}}

\newcommand{\tR}{{{\check R}}}
\newcommand{\tRuvx}{{\tR(u,v;x)}}
\newcommand{\tRbuvx}{{\tR^{\Box}(u,v;x)}}

\newcommand{\ve}{\varepsilon}
\newcommand{\rx}{{\rho_x}}
\newcommand{\ruvx}{{\rho_{u,v,x}}}\newcommand{\rvux}{{\rho_{v,u,x}}}
\newcommand{\Vx}{{V_x}}\newcommand{\Vqx}{{V_{qx}}}
\newcommand{\Vxi}{{V_x^{(i)}}}
\newcommand{\Vxichi}{{V_x^{(1)}}}\newcommand{\Vxni}{{V_x^{(2)}}}
\newcommand{\END}{{\rm End}}

\newcommand{\sis}{{s}}
\newcommand{\sig}{{\sigma}}

\newcommand{\lamb}{{\lambda}}

\newcommand{\Pxi}{{P_x^{(i)}}}
\newcommand{\Pxichi}{{P_x^{(1)}}}\newcommand{\Pxni}{{P_x^{(2)}}}
\newcommand{\chiu}{{\chi_u}}
\newcommand{\hnqtwo}{{H_n(q^2)}}
\newcommand{\Sn}{{S_n}}

\newcommand{\hi}{{h_i}}
\newcommand{\hs}{{h(\sigma)}}
\newcommand{\sigi}{{\sigma_i}}
\newcommand{\sigp}{{\sigma^\prime}}

\newcommand{\ep}{{e_+}}\newcommand{\ema}{{e_-}}
\newcommand{\epm}{{e_\pm}}

\newcommand{\Wxn}{{W_x^{(n)}}}\newcommand{\pixn}{{\pi_x^{(n)}}}
\newcommand{\tRiuvx}{{\tR_i(u,v;x)}}
\newcommand{\tRaxs}{{\tR(a;x|\sigma)}}

\newcommand{\rax}{{\rho_{a,x}}}

\newcommand{\tRp}{{{\check R}^\prime}}
\newcommand{\tRpp}{{{\check R}^{\prime\prime}}}

\newcommand{\ppm}{{p_\pm}}
\newcommand{\gamn}{{\gamma_n}} 
\newcommand{\Vpmx}{{V_{\pm,x}}} 
\newcommand{\dpmn}{{d_\pm(n)}}  
\newcommand{\ruxpmn}{{\rho_{u,x}^{\pm,(n)}}}
\newcommand{\tRpmnuvx}{{{\check R}^{\pm,(n)}(u,v;x)}}
\newcommand{\rpmnuv}{{\rho^{\pm,(n)}_{u,v}}}

\newcommand{\BrackenGouldZhang}{{4}}
\newcommand{\BazhanovShadrikov}{{8}}
\newcommand{\Cherednik}{{13}}
\newcommand{\CurtisReiner}{{16}}
\newcommand{\Drinfeld}{{2}}
\newcommand{\EtingofVarchenko}{{19}}
\newcommand{\Gade}{{9}}
\newcommand{\Gyoja}{{17}}
\newcommand{\IoharaKoga}{{7}}
\newcommand{\Jimboichi}{{3}}
\newcommand{\Jimbo}{{14}}
\newcommand{\Kac}{{5}}
\newcommand{\Kactwo}{{6}}
\newcommand{\KhoroshkinTolstoy}{{11}}
\newcommand{\KulishReshetikhinSklyanin}{{15}}
\newcommand{\PerkSchultz}{{10}}
\newcommand{\Sergeev}{{18}}
\newcommand{\Yamane}{{12}}
\newcommand{\YangGe}{{1}}

\newcommand{\lmmLmmone}{{1}\,\,}
\newcommand{\lmmLmmtwo}{{2}\,\,}
\newcommand{\lmmpropositionone}{{3}\,\,}
\newcommand{\thmtheoremone}{{1}\,\,}
\newcommand{\thmtheoremtwo}{{2}\,\,}

\begin{document}

 \begin{flushleft}
{\bf  {A central extension of $\Uqslninione$ 
and $R$-matrices with a new parameter}} 
\end{flushleft}

Hiroyuki Yamane
\par
Department of Pure and Applied Mathematics, Graduate School of Information Science
and Technology, Osaka University, Toyonaka 560-0043, Japan \par
e-mail:yamane@ist.osaka-u.ac.jp
\newline\newline
\begin{abstract} 
In this paper, using a quantum
superalgebra associated with the universal 
central extension of $\slninione$, we introduce
new $R$-matrices having an extra parameter $x$.
As $x\rightarrow 0$, they become
those associated with 
the symmetric and anti-symmetric tensor products of the copies of the
vector representation of $\Uqslninione$. 
\end{abstract}
(Received:
\newpage 
\section*{{\sf { I.
INTRODUCTION}}}
\quad\, The Yang-Baxter equation (YBE for short) has played important roles in study of
statistical mechanics, knot theory, conformal field theory 
etc.,$^\YangGe$ and many of its solutions
are associated with finite dimensional irreducible representations of quantum affine
algebras$^{\Drinfeld, \Jimboichi}$ and superalgebras.$^\BrackenGouldZhang$ We call the
solutions of  the YBE the {\it $R$-matrices.} \par 
 If a finite dimensional
simple Lie superalgebra is 
$A(m,n)$, $B(m,n)$, $C(n)$, $D(m,n)$, $F(4)$, $G(3)$
or $D(2,1;\al)$ ($\al\ne 0,\,-1$), it is called 
a {\it basic classical Lie superalgebra$^{\Kac, \Kactwo}$}  
(BCLS for short).
We first recall that 
$A(m,n)$ coincides with ${\rm sl}(m+1|n+1)$
if and only if $m\ne n$, and that ${\rm sl}(m+1|m+1)$
is a one dimensional central extension of $A(m,m)$.
Let $\G$ be a BCLS and $\ovG$ the universal central extension (UCE for
short) of $\G$. We also recall that
$\ovG=\G$ if
$\G\ne A(m,m)$ for any $m$, and that ${\overline {A(m,m)}}={\rm
sl}(m+1|m+1)$ ($m\geq 2$) and ${\overline {A(1,1)}}=\D$.
Here $\D$ is the Lie superalgebra called 
$D(2,1;-1)$.$^\IoharaKoga$ The
$\D$ is  a two (resp. three) dimensional central extension of ${\rm
sl}(2|2)$ (resp. $A(1,1)$).
The UCE of $\G\otimes\Ct$
is given by the affine version $\ovG^{(1)}=\ovG\otimes\Ct\oplus\C c$
of $\ovG$.$^\IoharaKoga$ Motivated by this fact, we direct our attention 
to the quantum superalgebra
$\UqDone$ (strictly speaking, $\wtU=\wtUqDone$)
in order to give new $R$-matrices $\tRuvx$ satisfying 
the (twisted) YBE:
\begin{eqnarray}
\lefteqn{
({\tR(v,w;x)}\otimes I)(I\otimes {\tR(u,w;q^nx)})
(\tRuvx\otimes I)}
\label{eqn:defYBE} 
\\
& &=(I\otimes {\tR(u,v;q^nx)})
({\tR(u,w;x)}\otimes
I)) (I\otimes {\tR(v,w;q^nx)}) 
\nonumber 
\end{eqnarray}
for some integer $n$, where $u$, $v$, $x\in\C$ are continuous parameters. 
This can be viewed as a quantum dynamical YBE
(see Appendix). 
The $R$-matrices we will give are such that as $x\rightarrow 0$,
they become the 
$\Uqslninione$ $R$-matrices$^{\BrackenGouldZhang, 
\BazhanovShadrikov{\rm -}\PerkSchultz
}$
associated  with the symmetric and anti-symmetric tensor products of
the copies of  the vector representation $\varphi$ of
$\Uqslninione$. 
One of our tools is a four dimensional irreducible representation
$\rx$ of $\wtU$ with the parameter $x$
such that 
$\rho_0=\varphi\circ p$, where 
$p:\wtU\rightarrow\Uqslninione$ is the natural epimorphism.   \par
The paper is organized as follows. 
In Section 1, we introduce $\wtU$ and $\rx$.
In Section 2, we give $\tRuvx$ associated with $\rx$.
In Section 3, we give all the $\tRuvx$'s mentioned above
using the fusion process.
\section*{{\sf {II. A CENTRAL EXTENSION OF $\Uqslninione$}}}
\quad\, Let ${\cal E}=\oplus_{i=0}^4\C\ve_i$ be the five dimensional
vector space. Define the symmetric bilinear form $(\,,\,)$
on ${\cal E}$ by $(\ve_0,\ve_0)=0$,
$(\ve_1,\ve_1)=(\ve_2,\ve_2)=1$,
$(\ve_3,\ve_3)=(\ve_4,\ve_4)=-1$ and $(\ve_i,\ve_j)=0$
($i\ne j$). Let $\al_0:=\ve_0-\ve_1+\ve_4$
and $\al_i:=\ve_i-\ve_{i+1}$ ($1\leq i\leq 3$).
Define the parity $p(\al_i)$ to be $(4-(\al_i,\al_i)^2)/4$. 
Then the Cartan matrix of $\Aoneoneone$ is
given by the $4\times 4$ matrix $(a_{ij})$,
where $a_{ij}=2(\al_i,\al_j)/((\al_i,\al_i)+2p(\al_i))$.
\par
Throughout this paper, we assume  
$q\in\C$ to be such that $q\ne 0$ and $q^r\ne 1$ 
for every positive
integer
$r$.  Let 
$\wtU=\wtUqDone$ be the associative $\C$-algebra presented by 
the generators $\sis$, $K_i^\pm$, $E_i$, $F_i$
($0\leq i\leq 3$) and the defining relations:
$$
\sis^2=1,\quad\sis K_i \sis=K_i,
\quad\sis E_i \sis=(-1)^{p(\al_i)}E_i,\quad\sis F_i
\sis=(-1)^{p(\al_i)}F_i,
$$
$$
K_iK_i^{-1}=1,\quad K_iK_j=K_jK_i,
$$
$$
K_iE_jK_i^{-1}=q^{(\al_i,\al_j)}E_j,\quad
K_iF_jK_i^{-1}=q^{-(\al_i,\al_j)}F_j,
$$
$$
[E_i,F_j]=\delta_{ij}{\frac {K_i-K_i^{-1}} {q-q^{-1}}}
\quad\mbox{if $(i,j)$ is neither $(2,0)$ nor $(0,2)$,}
$$
$$
K_2[E_2,F_0]\in Z(\wtU),\quad K_2^{-1}[E_0,F_2]\in Z(\wtU),
$$ where $[E_i,F_j]:=E_iF_j-(-1)^{p(\al_i)p(\al_j)}F_jE_i$
and $Z(\wtU)$ is the center of $\wtU$. We view $\wtU$
as the (non-$\Z_2$-graded) Hopf algebra with the comultiplication
$\Delta:\wtU\rightarrow\wtU\otimes\wtU$ satisfying:
$$
\Delta(\sis)=\sis\otimes\sis,\quad
\Delta(K_i)=K_i\otimes K_i,
$$
$$
\Delta(E_i)=E_i\otimes 1+
K_i\sis^{p(\al_i)}\otimes E_i
+\delta_{i0}(q-q^{-1})\sis [E_0,F_2]\otimes E_2,
$$
$$
\Delta(F_i)=F_i\otimes K_i^{-1}+
\sis^{p(\al_i)}\otimes F_i
-\delta_{i0}(q-q^{-1})F_2 \otimes [E_2,F_0].
$$ We do not give the antipode and the counit; we
do not need them. We define $\Delta^{(n-1)}:\wtU
\rightarrow \wtU^{\otimes n}$ by letting 
$\Delta^{(1)}=\Delta$ and $\Delta^{(m)}=
({\rm id}_\wtU\otimes \Delta^{(m-1)})\circ\Delta$ ($m\geq 2$).
\newline\newline{\it Remark:} (1) The above comultiplication
is not standard. Taking the twisting$^\KhoroshkinTolstoy$ for
the $\wtU$, we get the standard comultiplication
of a quantum superalgebra defined for a Dynkin diagram other than the one
associated with the Cartan matrix $(a_{ij})$ (see above); the ${\rm A}(1,1)^{(1)}$ has the
two Dynkin diagrams.
\par
(2) Let $\wtUp$ be the subalgebra of 
$\wtU$ generated by $K_i^\pm$, $E_i$, $F_i$. Then
$\wtU=\wtUp\oplus\wtUp\sis$.  There exists a nonzero ideal $J$ of $\wtUp$
such that
$\wtUp/J$ can be regarded as $\UqDone$. We can get generators of 
$J$ in the same way as in Ref.~\Yamane. By the same 
argument as in the proof of Theorem~8.4.3 of Ref.~\Yamane, 
we can get the natural epimorphism from $\UqDone$ to 
$U_q{\rm sl}(2|2)^{(1)}$.
\newline\par
Let $\Vx=\C^4$ be the four dimensional vector space,
where $x\in\C$ is a parameter. Put $\theta(i):= (1-(\ve_i,\ve_i))/2$.
Define the irreducible representation $\rx:\wtU\rightarrow\END(\Vx)$ by:
$$
\rx(\sis)=\sum_{j=1}^4(-1)^{\theta (j)}E_{jj},
\quad \rx(K_i)=\sum_{j=1}^4q^{(\al_i,\ve_j)}E_{jj},
$$
$$
\rx(E_0)=E_{41},\quad  \rx(E_1)=E_{12},
$$ 
$$  \rx(E_2)=E_{23}+x E_{41},
\quad  \rx(E_3)=E_{34},
$$
$$
\rx(F_0)=-E_{14}-xq^{-1}E_{32},
\quad  \rx(F_1)=E_{21},
$$ 
$$
\rx(F_2)=E_{32},
\quad  \rx(F_3)=-E_{43}.
$$ \par
Let $\wtUzero=\wtUqD$ be the subalgebra of $\wtU$ generated by
$\sis$, $K_i^\pm$, $E_i$, $F_i$ ($1\leq i\leq 3$).
Define the vector subspaces $\Vxi=V_{x,y}^{(i)}$ ($i=1,\,2$) of
$\Vx\otimes V_y$ by
\begin{eqnarray*}
\Vxichi&:=& {\C}(e_3\otimes e_3)\oplus{\C}(e_4\otimes e_4) \\
& &\oplus\bigoplus_{i<j}{\C}
(e_i\otimes e_j-(-1)^{\theta(i)\theta(j)}qe_j\otimes e_i \\
& &+\delta_{i1}\delta_{j2}(q^2ye_3\otimes e_4
+xe_4\otimes e_3))
\end{eqnarray*} and
$$
\Vxni := {\C}(e_1\otimes e_1)\oplus{\C}(e_2\otimes e_2)
\oplus\bigoplus_{i<j}{\C}
(e_i\otimes e_j+(-1)^{\theta(i)\theta(j)}q^{-1}e_j\otimes e_i).
$$

{\it Lemma 1:
$\Vxichi$ is an irreducible $\wtUzero$-module. Moreover
 $\Vx\otimes V_y$ is a completely reducible
$\wtUzero$-module
if and only if $y=qx$. If this is the case, 
$\Vxni$ is an irreducible $\wtUzero$-module
which is not isomorphic to $\Vxichi$;
in particular, 
$\Vx\otimes\Vqx$ has an irreducible $\wtUzero$-submodule
decomposition 
$\Vxichi\oplus\Vxni$.}
\par {\it Proof:} For each $1\leq i\leq 4$, the weight space including
$e_i\otimes e_i$ is one dimensional. 
Hence, if  $\Vx\otimes V_y$ is a completely reducible
$\wtUzero$-module,
there exists an irreducible $\wtUzero$-module
including $e_i\otimes e_i$. Using this fact,
we can check the lemma directly. \hfill $\Box$ 
\section*{{\sf {III. $R$-MATRIX FOR THE VECTOR REPRESENTATION}}}
\quad\, Define $\Pxi\in\END(\Vx\otimes\Vqx)$ ($i=1,\,2$)  
by $\Pxi(v)=\delta_{ij}v$ ($v\in V_x^{(j)}$). Put
\begin{equation}
\tRbuvx :=(q^2u-v)\Pxichi +(q^2v-u)\Pxni
, \label{eqn:defvR}
\end{equation} where $u,\,v\in\C$.  
Then:
\begin{eqnarray}
 \lefteqn{
\tRbuvx} \label{eqn:concreteR}\\
&=&(q^2v-u)\sum_{i=1}^2E_{ii}\otimes E_{ii}
+(q^2u-v)\sum_{i=3}^4E_{ii}\otimes E_{ii} \nonumber\\
& &+(q^2-1)\sum_{i<j}(vE_{ii}\otimes E_{jj}
+uE_{jj}\otimes E_{ii}) \nonumber\\
& &-q(u-v)\sum_{i\ne j}(-1)^{\theta(i)\theta(j)}E_{ij}\otimes E_{ji}
\nonumber\\ & &+x(q^2-1)(u-v)(qE_{31}\otimes E_{42}
-q^2E_{32}\otimes E_{41} \nonumber\\
& &-E_{41}\otimes E_{32}+ q E_{42}\otimes E_{31}).
\nonumber
\end{eqnarray}
\par For $u\in\Ctms$, define $\chiu\in{\rm Aut}(\wtU)$ by
$\chiu(\sis)=\sis$, $\chiu(K_i)=K_i$,  
$\chiu(E_i)=u^{-\delta_{i0}}E_i$ and
$\chiu(F_i)=u^{\delta_{i0}}F_i$. Put
$\ruvx:=(\rx\otimes\rho_{qx})\circ(\chiu\otimes\chi_v)\circ\Delta$.
Using (\ref{eqn:defvR}) and 
Lemma~\lmmLmmone,
we can directly 
check that: 
\begin{equation}
\tRbuvx\ruvx(X)=\rvux(X)\tRbuvx \label{eqn:RrrR}
\end{equation} for $X\in\wtU$. 

{\it {\bf Theorem 1:} The $\tRbuvx$ satisfies the YBE in the form of 
(\ref{eqn:defYBE}) with $n=1$. }
\par {\it Proof:} Let $E^\prime_i$, $F^\prime_i$, $H^\prime_i$
($0\leq i\leq 3$) be the Chevalley generators of 
${\rm sl}(2|2)^{(1)}$. Then there exists a 
representation ${\widehat \psi}_u:{\rm sl}(2|2)^{(1)}
\rightarrow\END(\C^4)$
sending $E^\prime_i$, $F^\prime_i$, $H^\prime_i$
to the limits of
$\rx\circ\chi_u(E_i)$, $\rx\circ\chi_u(F_i)$,
$(q-q^{-1})^{-1}\rx\circ\chi_u(K_i-K_i^{-1})$ as 
$(q,x)\rightarrow (1,0)$,  
respectively. Notice that 
$\psi:={\widehat \psi}_1{}_{|{\rm sl}(2|2)}$ 
is an irreducible representation of ${\rm sl}(2|2)$
and that there exist a highest root vector 
$E^\prime_{\al_1+\al_2+\al_3}$ and a lowest root vector
$E^\prime_{-(\al_1+\al_2+\al_3)}$ of ${\rm sl}(2|2)$
such that 
\begin{equation}  \label{eqn:upsi}
u\psi(E^\prime_{\al_1+\al_2+\al_3})
={\widehat \psi}_u(F_0),\quad
u^{-1}\psi(E^\prime_{-(\al_1+\al_2+\al_3)})
={\widehat \psi}_u(E_0).
\end{equation}
Then, using (\ref{eqn:RrrR}), together with the same argument used in
the proof of
Proposition~3 in Ref.~\Jimboichi, we get the theorem. 
\hfill $\Box$

\section*{{\sf {IV. $R$-MATRIX FOR THE (ANTI-)SYMMETRIC 
TENSORS}}} 
\quad\, Here we use a similar process to 
the fusion process.$ ^{\Cherednik-\KulishReshetikhinSklyanin
}$ To begin with, we recall some facts$^{\CurtisReiner, \Gyoja}$ about the
Hecke algebra
$\hnqtwo$ associated with the symmetric group $\Sn$; the $\hnqtwo$
is the associative $\C$-algebra presented by the generators 
$\hi$ ($1\leq i\leq n-1$) and the defining relations:
$(\hi -q^2)(\hi +1)=0$, $\hi h_{i+1}\hi=h_{i+1}\hi h_{i+1}$
and $\hi h_j=h_j\hi$ ($|i-j|\geq 2$).
We abbreviate $\hnqtwo$ to $H$.
We know that there exists a $\C$-basis
$\{\hs|\sig\in\Sn\}$ of $H$ such that $h(1)=1$, $h(\sigi)=\hi$ and
$h(\sigp\sig)=h(\sigp)\hs$ if $\ell(\sigp\sig)=\ell(\sigp)+\ell(\sig)$.
Here 
$\sigi$  is the simple transposition $(i,i+1)$
and $\ell(\sig)$ is the length of $\sig$ with
respect to 
$\sigi$'s.\par
Put
$$
\ep:=\sum_{\sig\in\Sn}\hs,\quad\ema:=\sum_{\sig\in\Sn}(-q^{-2})^{\ell(\sig)}\hs. 
$$ Then $\hi\ep=q^2\ep$, $\hi\ema=-\ema$ and

\begin{equation}
\epm^2=(\sum_{\sig\in\Sn}q^{\pm 2\ell(\sig)})\epm.
\label{eqn:epmsquare}
\end{equation} 

Now we treat $R$-matrices. Let $\Wxn :=\Vx\otimes\Vqx\otimes\cdots\otimes
V_{q^{n-1}x}$. Put: 
$$
\tRiuvx:=I^{\otimes i-1}\otimes\tR^\Box(u,v;q^{i-1}x)\otimes
I^{\otimes n-i+1}\in\END(\Wxn ).
$$ By 
Theorem~\thmtheoremone, we can define 
$\tRaxs\in\END(\Wxn )$, $a\in(\C^\times)^n$ and 
$\sig\in\Sn$, inductively by  
$$
\tR(a;x|1)=I^{\otimes n},\quad \tR(a;x|\sigi)=\tR_i(a_i,a_{i+1};x)
$$ and 
$$
\tR(a;x|\sigp\sig)=\tR(\sig [a];x|\sigp)\tRaxs
\quad\mbox{if $\ell(\sigp\sig)=\ell(\sigp)+\ell(\sig)$,}
$$ where $\sig[a]:=(a_{\sig^{-1}(1)},\ldots,a_{\sig^{-1}(n)})$.
By 
Theorem~\thmtheoremone
 and (\ref{eqn:defvR}), 
there exists a unique representation
$\pixn :H\rightarrow\END (\Wxn )$ such that 
$\tRiuvx=\pixn (u\hi-vq^2\hi^{-1})$. \par
Let $\ppm:=(1,q^{\mp 2},\ldots,q^{\mp 2(n-1)})\in\C^n$. 
Let $\gamn\in\Sn$ be such that 
$\gamn(i)=n-i+1$. 

{\it Lemma 2:  Let $u\in\C$. Then:
$$
\tR(u\ppm;x|\gamn)=u^{\ell(\gamn)}a_\pm(q)\pixn (\epm)
$$ for some $a_\pm(q)\in\Ctms$.}
\par
This can be checked directly; a similar formula
has been given in Section~5 in Ref.~\Jimbo. 
\newline\par
Let $\Vpmx:=\pixn(\epm)\Wxn$.
By (\ref{eqn:epmsquare}),  $\dpmn:=\dim\Vpmx$ does not depend on $q$ or
$x$. For $a\in(\Ctms)^n$, define the representation
$\rax :\wtU\rightarrow \END (\Wxn )$ by
$$
\rax:=(\rx\otimes\cdots\otimes\rho_{q^{n-1}x})\circ
(\chi_{a_1}\otimes\cdots\otimes\chi_{a_n})\circ
\Delta^{(n-1)}.
$$ \par By (\ref{eqn:RrrR}), we have:
\begin{equation}
\tRaxs\rax (X)=\rho_{\sig[a],x}(X)\tRaxs
\label{eqn:bigRrrR}
\end{equation} for $X\in\wtU$. By 
Lemma~\lmmLmmtwo
 and
(\ref{eqn:bigRrrR}), we may define the representation
$\ruxpmn:\wtU\rightarrow\END(\Vpmx)$ by
$\ruxpmn(X)=\rho_{\gamn[u\ppm],x}(X)_{|\Vpmx}$.
Notice that
\begin{equation}
\ruxpmn=\rho_{1,x}^{\pm,(n)}\circ\chiu .\label{eqn:ruxlxii}
\end{equation}

We
have a representation 
${\widehat \psi}^{\pm,(n)}_u:
{\rm sl}(2|2)^{(1)}\rightarrow \END (\C^\dpmn)$
sending $E^\prime_i$, $F^\prime_i$, $H^\prime_i$
to the limits of
$\ruxpmn(E_i)$, $\ruxpmn(F_i)$,
$(q-q^{-1})^{-1}\ruxpmn(K_i-K_i^{-1})$ as 
$(q,x)\rightarrow (1,0)$.
Define the representation $\psi^{\pm,(n)}$ 
of ${\rm sl}(2|2)$ to be
$({\widehat
\psi}^{\pm,(n)}_1) _{|{\rm sl}(2|2)}$. 
Then $\psi^{+,(n)}$ (resp. $\psi^{-,(n)}$)
is the $n$-fold symmetric (resp. anti-symmetric) tensor product
of the vector representation $\psi$ of ${\rm sl}(2|2)$.
By Ref.~\Sergeev, we have:

{\it Lemma 3: The $\psi^{\pm,(n)}$
is irreducible. Moreover $\dpmn\ne 0$. 
}

Define $\tau\in S_{2n}$ by $\tau (i)=i+n$, 
$\tau (n+i)=i$ ($1\leq i\leq n$).  
 For
$g,\,h\in\C^n$, let $g\cup
h:=(g_1,\ldots,g_n,h_1,\ldots,h_n)\in\C^{2n}$. 
Let $S_n$ be embedded into 
$S_{2n}$ in the natural way.
By 
Lemma~\lmmLmmtwo, we
have: 
\begin{eqnarray*}
\lefteqn{\tR(\gamn[u\ppm]\cup\gamn[v\ppm];x|\tau)
(\pixn(\epm)\otimes\pi_{q^nx}^{(n)}(\epm))} \\ 
&=& {\frac {(uv)^{-\ell (\gamn)}}
{a_\pm(q)^2}}\tR(\gamn[u\ppm]\cup\gamn[v\ppm];x|\tau)
\tR(u\ppm;x|\gamn)\tR(v\ppm;x|\tau\gamn\tau) \\
&=& {\frac {(uv)^{-\ell (\gamn)}}{a_\pm(q)^2}}
\tR(u\ppm\cup v\ppm;x|\gamn\tau\gamn) \\
&=& {\frac {(uv)^{-\ell (\gamn)}}{a_\pm(q)^2}}
\tR(v\ppm;x|\gamn)\tR(u\ppm;x|\tau\gamn\tau)\tR(u\ppm\cup
v\ppm;x|\tau) \\ &=& 
(\pixn(\epm)\otimes\pi_{q^nx}^{(n)}(\epm))
\tR(u\ppm\cup
v\ppm;x|\tau)\,.
\end{eqnarray*}
Hence we may put:
$$
\tRpmnuvx:=\tR(\gamn[u\ppm]\cup\gamn[v\ppm];x|\tau)
_{|\Vpmx\otimes
V_{\pm,q^nx}}
$$ $\in\END (
\Vpmx\otimes 
V_{\pm,q^nx}).
$

Let $\rpmnuv:=(\ruxpmn\otimes\rho^{\pm,(n)}_{v,q^nx})\circ\Delta$.
Notice that
\begin{equation} \label{eqn:rlrm}
\rpmnuv (X) =
(\rho_{\gamn[u\ppm]\cup\gamn[v\ppm],x}(X))
_{|\Vpmx\otimes V_{\pm,q^nx}}. 
\end{equation}
By (\ref{eqn:bigRrrR}) and (\ref{eqn:rlrm}), we have:
\begin{equation} \label{eqn:RrrRsononi}
\tRpmnuvx\rpmnuv (X)
=\rho^{\pm, (n)}_{v,u}
(X)\tRpmnuvx.  
\end{equation}  
for $X\in\wtU$.

{\it {\bf Theorem 2:} The $\tRpmnuvx$ satisfies the YBE in the form of (\ref{eqn:defYBE}).}
\par{\it Proof:} By (\ref{eqn:upsi}), we have
$u\psi^{\pm,(n)}(E^\prime_{\al_1+\al_2+\al_3})
={\widehat \psi}^{\pm,(n)}_u(F_0)$
and \newline $u^{-1}\psi^{\pm,(n)}(E^\prime_{-(\al_1+\al_2+\al_3)})$
$=$ ${\widehat \psi}^{\pm,(n)}_u(E_0)$.
Noting this fact and
using (\ref{eqn:ruxlxii}), (\ref{eqn:RrrRsononi}) 
and 
Lemma~\lmmpropositionone, together
 with the same argument
as in the proof of Proposition 
3  in Ref.~\Jimboichi,  we have
the theorem. \hfill $\Box$
\newline\newline
\section*{{{\sf ACNOWLEDGMENTS}}}
\quad\, The author thanks E.~Date, M.~Okado and Y.~Koga
for valuable comments.
He also thanks Y.Z.~Zhang for telling him about
face-type dynamical $R$-matirces.
\section*{{\sf {APPENDIX: A QUANTUM DYNAMICAL $R$-MATRIX}}}  
\quad\, Here we show that the $\tRpmnuvx$ can be viewed as a dynamical 
$R$-matrix. Let ${\frak{h}}$ be a finite dimensional
commutative Lie algebra. 
Let $V$ be a finite dimensional diagonalizable ${\frak{h}}$-module, i.e.,
$V=\oplus_{\mu\in{\frak{h}}^*}V_\mu$, where $V_\mu:=\{v| h.v=\mu(h)v\}$.
We say that a (meromorphic) function 
$\tRp:{\C}^2\times{{\frak{h}}}^*\rightarrow\END(V\otimes V)$ 
is a {\it quantum dynamical $R$-matrix} if it satisfies the
quantum dynamical YBE (see
Ref.~\EtingofVarchenko\, for example):
\begin{eqnarray*}
\lefteqn{(\tRp(v,w,\lamb)\otimes I)\tRp_{23}(u,w,\lamb
-h^{(1)})(\tRp(u,v,\lamb)\otimes I)}
\\
& &=\tRp_{23}(u,v,\lamb
-h^{(1)})
({\tRp(u,w,\lamb)}\otimes
I) \tRp_{23}(v,w,\lamb
-h^{(1)}), 
\end{eqnarray*}
where $\tRp_{23}(u,v,\lamb
-h^{(1)})\in\END(V^{\otimes 3})$ is defined by 
$$\tRp_{23}(u,v,\lamb
-h^{(1)})_{|V_\mu\otimes V\otimes V}
=(I\otimes\tRp(u,v,\lamb
-\mu))_{|V_\mu\otimes V\otimes V}.
$$ 
\par
Let ${{\frak{h}}}^{\prime\prime}
={\C}$ and let ${{\frak{h}}}^{\prime\prime}$ act on 
$\C^\dpmn$ by $z.v=-nzv$. Let $a\in\C$ be such that
$e^a=q$.
Define
$\tRpp:{\C}^2\times ({\frak{h}}^{\prime\prime})^*
\rightarrow\END(\C^\dpmn\otimes\C^\dpmn)$
by $\tRpp(u,v,\lamb)=\tR^{\pm,(n)}(u,v;e^{a\lamb(1)})$. 
By 
Theorem~\thmtheoremtwo, $\tRpp$ is a quantum dynamical
$R$-matrix.  
\newline\newline\newline
$^\YangGe$C.~N.~Yang and M.~L.~Ge(eds),
\textit{Braid Groups, Knot theory and Statistical Mechanics}, (World 
Scientific, Singapore, 1989). \newline
$^\Drinfeld$V.~G.~Drinfeld,
\textit{ICM Proceedings,} (Berkeley 1986), p.~798.
\newline
$^\Jimboichi$M.~Jimbo, 
{\it Commun. Math. Phys.}
{\bf 102}, 537 (1986).
\newline
$^\BrackenGouldZhang$A.~J.~Bracken, M.~D.~Gould and R.~B.~Zhang,,
{\it Mod. Phys. Lett. A} {\bf 5}, 831 (1989).
\newline
{$^\Kac$V.~G.~Kac, 
{\it Adv. in. Math.}  
{\bf 26}, 8 (1977).
\newline
{$^\Kactwo$V.~G.~Kac, 
{\it Lect. Note in Math.}  
{\bf 676}, 597 (1978).
\newline
$^\IoharaKoga$K.~Iohara and Y.~Koga, 
{\it Comment. Math. Helv.}  
{\bf 76}, 110 (2001)
\newline
$^\BazhanovShadrikov$V.~V.~Bazhanov and A.~G.~Shadrikov,
{\it Theoret. Math. Phys.} {\bf 73}, 1302 (1987).
\newline
$^\Gade$R.~M.~Gade, 
{\it J. Phys. A.}
{\bf 31}, 4909 (1998). 
\newline
$^\PerkSchultz$J.~H.~H.~Perk, and C.~L.~Schultz,
{\it Non-linear Integrable systems--Classical theory and 
quantum theory} ed, by M.~Jimbo and
T.~Miwa, (World Scientific, Singapore, 1983), p.~135.
\newline
$^\KhoroshkinTolstoy$S.~M.~Khoroshkin and V.~N.~Tolstoy,
{\it Twisting of quantum (super-)algebras,}
{\it in Generalized symmetries in physics}
(Claustal, 1993), (World Sci. Publishing, River Edge, NJ,
1994) p.~42.
\newline 
$^\Yamane$H.~Yamane,
{\it Publ. RIMS Kyoto Univ.}  
{\bf 35}, 321 (1999),
(Errata) {\it Publ. RIMS Kyoto Univ.}  
{\bf 37}, 615 (2002).
\newline
$^\Cherednik$I.~V.~Cherednik, 
 {\it Sov. Math. Dokl.} 
{\bf 33}, 507 (1986).
\newline
$^\Jimbo$M.~Jimbo, 
{\it Lett. Math. Phys.}
{\bf 11}, 247 (1986).
\newline
$^\KulishReshetikhinSklyanin$P.~P.~Kulish, N.~Yu~Reshetikhin and E.~K.~Sklyanin,
{\it Lett. Math. Phys.}  
{\bf 5}, 393 (1981).
\newline
$^\CurtisReiner$C.~W.~Curtis and I.~Reiner, 
{\it Methods of Representation Theory 
with Applications to Finite Groups and Orders, Volume II},
(A Wiley International Publication, New York. 1987).
\newline
$^\Gyoja$A.~Gyoja, 
{\it Osaka J. Math.}
{\bf 23}, 841 (1986).
\newline
$^\Sergeev$A.~N.~Sergeev,
{\it Math. USSR-Sb.}  
{\bf 51}, 419 (1985).
\newline
$^\EtingofVarchenko$P.~Etingof and A.~Varchenko,
{\it Commun. Math. Phys.} {\bf 196}, 591 (1998). 
\end{document}